\documentclass[preprint,times]{elsarticle}
\usepackage{amsmath}
\journal{Journal of Computational Physics}

\renewcommand{\dot}[2]{\langle #1|#2\rangle}
\newcommand{\g}[2]{\langle #1|#2\rangle_\text{G}}

\def\Abf{\mbox{$\mathbf{A}$}}

\def\ubf{\mbox{$\mathbf{u}$}}

\def\vbf{\mbox{$\mathbf{v}$}}
\def\Sbf{\mbox{$\mathbf{S}$}}
\def\Dbf{\mbox{$\mathbf{D}$}}
\def\Mbf{\mbox{$\mathbf{M}$}}

\begin{document}

\begin{frontmatter}

\title{Short note on the mass matrix for Gauss-Lobatto grid points}

\author{Saul A. Teukolsky}
\ead{saul@astro.cornell.edu}

\address{Departments of Physics and Astronomy, Space
Sciences Building, Cornell University,
Ithaca, NY 14853, United States}

\begin{abstract}
The mass matrix for Gauss-Lobatto grid points is usually approximated
by Gauss-Lobatto quadrature because this leads to a diagonal matrix
that is easy to invert. The exact mass matrix and its inverse are full.
We show that the exact mass matrix \emph{and} its inverse differ from the
approximate diagonal ones by a simple rank-1 update (outer product).
They can thus be applied to an arbitrary vector in $O(N)$ operations instead
of $O(N^2)$.
\end{abstract}


\end{frontmatter}

\section{Motivation}
\label{sec:motiv}

With the increased emphasis on higher-order methods for solving partial
differential equations, methods that divide the domain into subdomains
and represent the solution as an
expansion in basis functions have become more and more important.
These include spectral element methods (penalty-based or continuous)
and discontinuous Galerkin methods.
To handle nonlinearities, collocation schemes are often the method of choice.
In such methods, the expansion coefficients are replaced by function
values at specially
chosen grid points as the fundamental unknowns.
In one dimension, the grid
points are universally chosen to be the Gaussian quadrature points
corresponding to the basis functions. This connects the expansion coefficients
in spectral space to the function values in physical space by a discrete
transform and leads to rapidly convergent and stable methods for smooth
solutions. 

In two and three dimensions, if the subdomains can be mapped to squares
or cubes, then basis functions that are tensor products of one-dimensional
basis functions are almost always used because of the resulting simplification
of element-wise operations. Unless the problem requires the flexibility of
grids constructed using triangles or tetrahedra, this approach is
again almost universal. The key result of this note applies to
any one-dimensional set of grid points that define a Gaussian quadrature or
are part of a tensor product of such grid points. It does not apply
to typical basis sets for triangles, where the quadrature rule and the choice
of grid points are not directly connected.

For many problems, the simplest formulation uses Gauss-Lobatto
collocation points since having grid points on the boundaries
makes it easy to impose boundary conditions. In such a formulation,
the exact mass matrix and its inverse are full.
Thus it is natural to approximate the mass matrix by Gauss-Lobatto quadrature,
which leads to a diagonal matrix that is easy to invert.
By contrast, using Gauss collocation points with ordinary Gauss
quadrature gives the exact mass matrix, which is diagonal. This makes
the comparison between the two choices tricky. On the one hand,
Gauss-Lobatto avoids interpolation from the interior points to the
boundaries, but on the other hand it may require more collocation points
to achieve the same accuracy as using Gauss points if you use
the approximate mass matrix for efficiency.
This point is discussed further in \S~\ref{sec:effic}.

We show that there is a simple expression for the exact mass matrix
and its inverse for Gauss-Lobatto collocation. Multiplying a vector
by one of these expressions can be done in $O(N)$ operations, just as for
a diagonal matrix. This suggests that efficiency versus accuracy
results for implementations of spectral methods should be reconsidered.
Of course, for large values of $N$ the spectral convergence of
Gaussian quadrature is likely to make the difference between the exact and
approximate mass matrices irrelevant. However, for small or moderate
$N$ the situation is not clear.

\section{Spectral Approximation}
This section summarizes some standard material
\cite{hesthaven2007,canuto2006,boyd2001,fornberg1996}
on spectral approximations in order to derive the key result in
the next section.

Consider approximations of functions
by expansions in orthogonal
polynomials:
\begin{equation}
u(x)=\sum_{k=0}^N b_k p_k(x)
\label{eq:modes}
\end{equation}
where
\begin{equation}
\int_{-1}^1 p_j(x) p_k(x)W(x)\,dx = h_k \delta_{jk}
\label{eq:orthog}
\end{equation}
The associated inner product is
\begin{equation}
\dot{u}{v}\equiv\int_{-1}^1 u(x)v(x)W(x)\,dx
\label{eq:inner}
\end{equation}
For simplicity, we will take the weight function
$W(x)=1$, in which case the basis functions are Legendre polynomials.
However, almost everything in this note goes through for
other systems of orthogonal polynomials.

The set of orthogonal polynomials determines a Gaussian quadrature
formula with weights $w_j$ and grid points $x_j$:
\begin{equation}
\int_{-1}^1 f(x)\,dx \approx \sum_{j=0}^N w_j f(x_j)
\label{eq:gauss}
\end{equation}
The Gauss-Lobatto version of this quadrature arranges for the
endpoints of the interval to be included in the set $x_j$.
Having collocation points on the boundary can make the
application of boundary conditions easier.
The quadrature \eqref{eq:gauss} is exact for polynomials of degree
no more than
$2N+1$ for the Gauss case and $2N-1$ for the Gauss-Lobatto case.
Use the quadrature to define the discrete inner product as the analog
of \eqref{eq:inner}:
\begin{equation}
\g{u}{v}=\sum_{j=0}^N w_j u(x_j) v(x_j)
\label{eq:gaussinner}
\end{equation}
The continuous and discrete inner products are the same if the
product $uv$ is a polynomial of degree no more than $2N+1$ (Gauss) or
$2N-1$ (Gauss-Lobatto).

Equation \eqref{eq:modes} is called a modal expansion.
In collocation methods, instead of regarding the $N+1$ modal coefficients
$b_k$ as fundamental, we choose a set of $N+1$ collocation points
$x_j$. Typically these are the Gauss or Gauss-Lobatto points
associated with the orthogonal polynomials.
The corresponding nodal expansion is
\begin{equation}
u(x)=\sum_{j=0}^N u_j \ell_j(x)
\label{eq:nodes}
\end{equation}
where $u_j\equiv u(x_j)$. The basis functions $\ell_j(x)$ are called
cardinal functions and are simply the Lagrange interpolating polynomials
based on the grid points $x_j$, with $\ell_j(x_i)=\delta_{ij}$:
\begin{equation}
\ell_j(x)=\prod_{\substack{i=0\\i\neq j}}^N \frac{x-x_i}{x_j-x_i}
\end{equation}

The nodal expansion \eqref{eq:nodes} is just an approximation of
a continuous function $u(x)$ by its interpolating polynomial, so that
$u(x_i)=u_i$. Note that in the discrete inner product \eqref{eq:gaussinner}
of any two continuous functions, we may replace $u$, say, by its
interpolating polynomial, since only the values $u_j$ contribute to the sum.
Thus with collocation methods we don't distinguish between
a function and its expansion when using discrete inner products.

Since the discrete and continuous inner products are the same for
polynomial integrands up to degree $2N-1$, the $p_k$'s satisfy the
discrete orthogonality condition
\begin{equation}
\g{p_j}{p_k}=\gamma_k \delta_{jk}
\end{equation}
where
\begin{equation}
\gamma_k\equiv\g{p_k}{p_k}=\sum_{j=0}^N w_k p_k(x_j)^2
\label{eq:gammak}
\end{equation}
For Gauss-Legendre quadrature,
\begin{equation}
\gamma_k = h_k = 2/(2k+1)\quad\text{(Gauss-Legendre)}
\end{equation}
where $h_k$ is the normalization defined in \eqref{eq:orthog}:
\begin{equation}
h_k\equiv\dot{p_k}{p_k}
\label{eq:hk}
\end{equation}
This is because the integration in \eqref{eq:gammak} is
exact since the degree of the polynomial in the integrand in \eqref{eq:hk}
is $2N$.

For Gauss-Legendre-Lobatto (GLL) points, by contrast,
\begin{equation}
\gamma_k=\begin{cases}
2/(2k+1), & 0\leq k < N\\
2/N, & k=N
\end{cases}
\end{equation}
In this case, the degree of exactness is only $2N-1$, and so $\gamma_N\neq
h_N$. This simple fact is at the root of the ``difficulties'' of
using GLL points.

In a collocation method, we regard the expansion \eqref{eq:modes}
not as an independent alternative, but as the corresponding expansion of
the interpolating polynomial \eqref{eq:nodes}. Evaluating the
expression \eqref{eq:modes}  at
the grid points gives the relation
\begin{equation}
u_i=\sum_{k=0}^N b_k p_k(x_i)
\label{eq:s2p}
\end{equation}
This can be regarded as a discrete transform from spectral space,
characterized by the representation $b_k$,
to physical space, characterized by the
$u_i$. To find the inverse transform, consider
\begin{equation}
\g{u}{p_k}=
\sum_{j=0}^N b_j \g{p_j}{p_k} = \sum_{j=0}^N b_j \gamma_k \delta_{jk}
= b_k \gamma_k
\end{equation}
Thus
\begin{equation}
b_k=\frac{1}{\gamma_k}\g{u}{p_k}=\frac{1}{\gamma_k}\sum_{j=0}^N w_j
p_k(x_j) u_j
\label{eq:p2s}
\end{equation}
This is the transform from physical to spectral space.


The cardinal functions $\ell_j(x)$ are polynomials of degree $N$ and
so they can be expanded as
\begin{equation}
\ell_j(x)=\sum_{k=0}^N a_k p_k(x)
\label{eq:lexp}
\end{equation}
where by \eqref{eq:p2s}
\begin{equation}
a_k=\frac{1}{\gamma_k}\sum_{i=0}^N w_ip_k(x_i) \ell_j(x_i)
=\frac{1}{\gamma_k}\sum_{i=0}^N w_ip_k(x_i) \delta_{ij}
=\frac{1}{\gamma_k} w_j p_k(x_j)
\end{equation}
Substituting this in equation \eqref{eq:lexp} gives
\begin{equation}
\ell_j(x)=w_j\sum_{k=0}^N\frac{1}{\gamma_k}p_k(x_j)p_k(x)
\label{eq:card}
\end{equation}
This expansion for the cardinal functions will be extremely useful
in what follows.

\section{Exact Expressions for the Mass Matrix and Its Inverse}
\label{app:a}
\subsection{The Mass Matrix}
The mass matrix is defined as
\begin{equation}
M_{ij}=\int_{-1}^1 \ell_i(x)\ell_j(x)\,dx=\dot{\ell_i}{\ell_j}
\label{eq:mass}
\end{equation}
Here we have taken the range of $x$ to be $[-1,1]$.
The derivation below goes through even when a weight function $W(x)\neq 1$
is included in \eqref{eq:mass}.

Evaluating the mass matrix by Gaussian quadrature gives
a diagonal matrix:
\begin{equation}
M_{ij}=\sum_{k=0}^N w_k \ell_i(x_k)\ell_j(x_k)=
\sum_{k=0}^N w_k \delta_{ik}\delta_{jk}=w_i\delta_{ij}
\label{eq:massgauss}
\end{equation}
This expression is exact for Gaussian quadrature, but not for
the Gauss-Lobatto case because the integrand is of degree $2N$.

Let's derive an exact expression for the Lobatto case.
Substituting expression \eqref{eq:card} for the cardinal functions
gives
\begin{align}
M_{ij}
&=\sum_{k=0}^N\sum_{l=0}^N w_i w_j \frac{1}{\gamma_k\gamma_l}
p_k(x_i)p_l(x_j) \int_{-1}^1 p_k(x)p_l(x)\,dx\notag\\
&=\sum_{k=0}^N\sum_{l=0}^N w_i w_j \frac{1}{\gamma_k\gamma_l}
p_k(x_i)p_l(x_j)\delta_{kl}h_k\notag\\
&=\sum_{k=0}^N w_i w_j \frac{h_k}{\gamma_k^2} p_k(x_i)p_k(x_j)\notag\\
&=\sum_{k=0}^N w_i w_j \frac{1}{\gamma_k} p_k(x_i)p_k(x_j)
+\left(\frac{h_N}{\gamma_N^2}-\frac{1}{\gamma_N}\right)
w_i w_j p_N(x_i)p_N(x_j)\label{eq:massinterm}\\
&=w_i\ell_j(x_i)+\left(\frac{h_N}{\gamma_N^2}-\frac{1}{\gamma_N}\right)
w_i w_j p_N(x_i)p_N(x_j)\notag\\
&=w_i\delta_{ij}+\alpha w_i w_j p_N(x_i)p_N(x_j)
\label{eq:mass2}
\end{align}
where we have defined
\begin{equation}
\alpha=\frac{h_N-\gamma_N}{\gamma_N^2}
\end{equation}
Equation \eqref{eq:massinterm} follows from the previous line
because $\gamma_k=h_k$ for $k<N$.

Equation \eqref{eq:mass2} reduces to equation \eqref{eq:massgauss} if
$\gamma_N=h_N$, as for Gauss points. But we see that for the Lobatto case,
where it is convenient to use the diagonal expression \eqref{eq:massgauss}
because it is easy to invert, we introduce an error because of our
``quadrature crime.'' Since the error in applying the
approximate mass matrix to
a vector converges away spectrally
fast for smooth problems
as we increase $N$, it is customary to ignore this error because of the other
benefits of Lobatto points.  However, there is no need to do this:
the extra term in equation \eqref{eq:mass2} is proportional to
the outer product of the vector $w_i p_N(x_i)$ with itself. This means
that in applying the mass matrix to a vector in a matrix-vector multiply,
the extra term can be computed as a dot product of $w_j p_N(x_j)$ with the
vector and then a scaling of the vector $\alpha w_i p_N(x_i)$ by the dot
product.
The operation count is $O(N)$, the same as from the diagonal term
$w_i \delta_{ij}$.

More importantly, the inverse of the mass matrix is equally simple,
as we now show.

\subsection{Inverse Mass Matrix}
The inverse of the mass matrix follows from the Sherman-Morrison formula:
\begin{equation}
(\Abf + \ubf\otimes\vbf)^{-1}
= \Abf^{-1} -
     \frac{(\Abf^{-1}\cdot\ubf)\otimes(\vbf\cdot\Abf^{-1})}{1+
     \vbf\cdot\Abf^{-1}\cdot\ubf}
\end{equation}
In our case,
\begin{equation}
\Abf=\text{diag}(w_i),\quad u_i=w_i p_N(x_i),\quad v_i = \alpha u_i
\end{equation}
We find
\begin{equation}
M^{-1}_{ij}=\frac{1}{w_i}\delta_{ij}+\beta p_N(x_i)p_N(x_j),\qquad
\beta\equiv -\frac{h_N-\gamma_N}{\gamma_N h_N}
\label{eq:massinv}
\end{equation}
The simple form of the extra off-diagonal term in equation \eqref{eq:massinv}
makes it easy to use the exact inverse in applications. Once again, applying
the inverse matrix to a vector is an $O(N)$ operation.

\section{Applications}

\subsection{The Differentiation Matrix}
\label{sec:diffmatrix}
As a trivial application, consider the differentiation matrix
that appears when solving partial differential equations:
\begin{equation}
D_{ij}= \ell_j'(x_i)
\end{equation}
where a prime denotes a derivative. The differentiation matrix typically
appears via the stiffness matrix $\Sbf$:
\begin{equation}
\Dbf=\Mbf^{-1}\cdot\Sbf
\label{eq:diffmat}
\end{equation}
where
\begin{equation}
S_{jk}=\int_{-1}^1 \ell_j(x)\ell_k'(x)\,dx = \dot{\ell_j}{\ell_k'}
=\g{\ell_j}{\ell_k'}
\label{eq:sdef}
\end{equation}
The last equality follows since the degree of the polynomial in the integrand
of \eqref{eq:sdef} is $2N-1$. So carrying out the quadrature gives the
exact result
\begin{equation}
S_{jk}=\sum_m w_m \ell_j(x_m)\ell_k'(x_m) = \sum_m w_m \delta_{mj}\ell_k'(x_m)
= w_j \ell_k'(x_j)
\end{equation}

It is well known (e.g., \cite{carpenter1996})
that if the approximate mass matrix \eqref{eq:massgauss}
is used in \eqref{eq:diffmat}, one gets the exact result for the
differentiation matrix:
\begin{equation}
\sum_j\big(\Mbf^{-1}_{\rm GLL}\big)_{ij}S_{jk}=
\sum_j \frac{1}{w_i}\delta_{ij} w_j \ell_k'(x_j) = \ell_k'(x_i) =D_{ik}
\end{equation}
But why exactly do we get the right answer without using the exact
mass matrix? One way of seeing this is to show explicitly that the
``extra'' terms in \eqref{eq:massinv} give no additional contribution:
\begin{equation}
\sum_j\beta p_N(x_i)p_N(x_j)S_{jk} \propto \sum_jp_N(x_j)w_j \ell_k'(x_j)
= \g{p_N}{\ell_k'} = 0
\end{equation}
Here the quadrature gives zero by orthogonality because the degree of
$\ell_k'$ is less than the degree of $p_N$.

\subsection{Projection in hp-refinement}
An advantage of methods like the DG method is that it is relatively
straightforward to implement adaptive mesh refinement, including
full $hp$-refinement. With refinement, there are two methods for communicating
fluxes across subdomain faces: interpolation and projection. Interpolation
is simpler, but for marginally resolved problems the inherent aliasing
can lead to an instability. Moreover, interpolation does not guarantee
conservation and so can be less robust than projection, especially for problems
with shocks.

A convenient way to implement projection is with mortars, auxiliary
slices inserted at boundary interfaces. A full
discussion with explicit formulas is given in \cite{kopriva1996,kopriva2002}.
We note here that when projecting the solution from the subdomain to the
mortar to be able to compute the flux, one gets the exact projection
matrix using Gauss-Lobatto quadrature even when using the approximate
mass matrix.
The proof of this result is similar to that of \S\ref{sec:diffmatrix}:
the extra outer product terms give no contribution. The resulting expression
then shows that the result is the \emph{same} as using
interpolation. It is only when transferring the flux back from the mortar
to the subdomain that there is a difference between projection
and interpolation, and only when the polynomial degree of the subdomain
is less than that of the mortar.
This observation can be used to greatly simplify the implementation
of projection for DG fluxes as given, for example, in
\cite{kopriva1996,kopriva2002}.
In retrospect, it is ``obvious'' that
interpolation from a coarse grid to a finer one introduces
no aliasing, and so projection
should be the same as interpolation, but this fact has not been used before
in the literature on $hp$-refinement, to the best of my knowledge.

\subsection{Efficiency of finite element methods}
\label{sec:effic}
The question of whether to use Gauss points or
Gauss-Lobatto points in spectral collocation methods is not always
clear-cut, especially for small or moderate numbers of grid points.
Gauss points typically require interpolation to impose boundary
conditions, but their higher degree of exactness may allow a smaller
number of points to be used for a given accuracy.
There have been a number of studies of this
question \cite{kopriva2010,bassi2013}.
Kopriva and Gassner \cite{kopriva2010} concluded that for a simple
linear wave equation, the two approximations have comparable efficiency,
but for a nonlinear steady-state example Gauss approximation was faster
for a desired error.
Bassi et al. \cite{bassi2013} concluded that Gauss
nodes have ``a clear advantage'' for steady-state problems, and
ascribed this to under-integration with Gauss-Lobatto nodes.
As already mentioned in \S\ref{sec:motiv}, this question should be
re-examined in the light of the result of this paper, which allows the
exact mass matrix or its inverse to be used efficiently in the Lobatto
case.

\subsection{Dispersion and dissipation}
For wave propagation problems, dispersion and dissipation errors are
important properties of any numerical scheme. Gassner and
Kopriva \cite{gassner2011} showed that the error introduced by
using the approximate Gauss-Lobatto mass matrix can be interpreted
as a modal filter applied to the highest polynomial mode (since $\gamma_N \neq
h_N$). This filtering greatly increases the dispersion
and dissipation errors compared to the Gauss case. It would be
worthwhile to re-examine this question with the exact mass matrix.

\subsection{Roundoff errors}
It has been noted in \cite{gottlieb2009} that different
ways of computing the terms in a spectral element or DG method
can affect roundoff errors as $N$ increases. It may be worth examining
whether the numerical behavior is affected by the different ways of
computing the mass matrix.

\section*{Acknowledgments}
I thank Jan Hesthaven for several helpful comments, including pointing
out that interpolation from a coarse grid to a finer one introduces no
aliasing.
This work was supported in part by
NSF Grants PHY-1306125 and AST-1333129 at Cornell University, and
by a grant from the Sherman Fairchild Foundation.

\section*{References}

\end{document}